\documentclass[a4paper,10pt]{article}
\usepackage{amsmath, amsthm, amssymb, amsfonts, epsfig, epstopdf, titling, url, array, enumerate, cite, mathrsfs, upgreek, authblk, scrextend, soul, color, graphicx, float, tikz}
\usepackage[bindingoffset=0.2in,left=1in,right=1in,top=1in,bottom=1in,footskip=.25in]{geometry}

 \usepackage[bookmarksnumbered, colorlinks, plainpages]{hyperref}
  \hypersetup{colorlinks=true,linkcolor=pink, anchorcolor=green, citecolor=cyan, urlcolor=orange, filecolor=magenta, pdftoolbar=true}

\theoremstyle{plain}
\newtheorem{thm}{Theorem}[section]
\providecommand{\keywords}[1]{\begin{addmargin}[28pt]{28pt}\noindent\textbf{Keywords:} #1 \end{addmargin}}
\newtheorem{lma}[thm]{Lemma}

\newtheorem{ppn}[thm]{Proposition}
\theoremstyle{definition}
\newtheorem{dfn}[thm]{Definition}
\newtheorem{eg}[thm]{Example}
\newtheorem{rem}[thm]{\textit{Remark}}
\providecommand{\ams}[1]{\begin{addmargin}[28pt]{28pt}\noindent\textbf{Mathematics Subject Classification:} #1\end{addmargin}}

\title{Exploring fixed point results  in fuzzy $\mathscr{F}$-metric spaces with  an   application  to satellite web coupling problem }
\author[1]{Dipti Barman}
\author[2]{Abhishikta Das}
\author[3]{T. Bag\thanks{Corresponding author}}
\affil{Department of Mathematics, Siksha-Bhavana,	\authorcr	Visva-Bharati, Santiniketan-731235, Birbhum, West-Bengal, India 
	\authorcr 	E-mail Id:  diptibarmanhmt@gmail.com$ ^1 $, abhishikta.math@gmail.com$^2 $, tarapadavb@gmail.com$ ^3 $
	\authorcr Orcid Id: 0009-0006-1146-4566$^1 $, 0000-0002-2860-424X$^2 $,    0000-0002-8834-7097$^3 $	}

\date{}

\begin{document}
	
 \maketitle

 \pagestyle{myheadings}
  	\markright{\footnotesize \it  Exploring fixed point results  in fuzzy $\mathscr{F}$-metric spaces with  an   application  to satellite web coupling problem}

\begin{abstract}
\noindent
This article explores several fundamental aspects of fuzzy $\mathscr{F}$-metric spaces and their applications in mathematical analysis. We  investigate  some essential  properties concerning compactness and total boundedness in  fuzzy $\mathscr{F}$-metric spaces. Within  this framework, we present a  fixed point theorem and demonstrate its utility by applying it to the satellite web coupling problem. To support the theoretical findings, illustrative examples and a graphical representation of the contraction condition are also provided. 
\end{abstract}

	\ams{46S40, 54E35.}
    
\section{Introduction}
In 1965, L. A. Zadeh \cite{5} made a   breakthrough contribution by introducing the theory of   fuzzy set. This   idea led to the development of fuzzy logic, fuzzy relations  and also to countless extensions such as fuzzy topology and fuzzy algebra. In this progress, introduction of fuzzy metric spaces was a natural progression that combines the principles of  classical metric spaces with fuzzy set theory. The concept of fuzzy metric spaces was first introduced by Kramosil and Michalek \cite{1}  in 1975. Later   in 1994,  George and Veeramani  \cite{2} refined this framework   incorporating continuous $t$-norms to make  the induced topology Hausdorff.   Day by day,  several generalizations of fuzzy metric spaces have been proposed viz. intuitionistic fuzzy metric space \cite{6}, $n$-dimensional fuzzy metric space \cite{7}, probabilistic fuzzy metric space \cite{8} and fuzzy 2-metric space \cite{9}. With this process of generalization, extension of fixed point theory has been made via fuzzy set theory(for references see \cite{11,12, AM-1,AM-2}).    
These advancements   have found applications in diverse areas such as   uncertainty model,    image processing, decision-making, optimization  and dynamic systems. 

 These successive generalizations have motivate  further extensions in the development of fuzzy  metric spaces. As a result,  Das et al. \cite{4} recently introduced   fuzzy $\mathscr{F}$-metric spaces   which seek to explore new theoretical properties and practical applications.   Such an extension  involves a function $\textit{f} : [ 0, 1 ] \to [ 0, 1 ] $ belonging to  a specific class of functions  to relax the axioms of   traditional fuzzy metric spaces. Incorporation of such function $\textit{f}  $ in this framework allows for greater flexibility in modeling and analysis. This generalized structure not only enlarge the    fuzzyness of distance functions    but also provides new ways in diverse directions for exploring fixed-point theorems and their applications.

Continuing this line of research, the present article focuses on investigating the fundamental properties of fuzzy $\mathscr{F}$-metric spaces and  the extension of   the development of fuzzy fixed-point results.  We establish few results concerning compactness, closed and bounded sets  and totallly boundedness. Additionally, we present new fixed point theorems under generalized contraction conditions, with applications to problems such as the existence and uniqueness of solutions to nonlinear ordinary differential equations. A geometric illustration is provided to support one of the key examples, highlighting the interplay between various parameters. Through appropriate examples, the necessity and sufficiency of the proposed contraction mappings are thoroughly justified.  
 
The structure of this paper is as follows.   Section 2 consists of   necessary preliminary results. In Section 3 we prove several basic  properties of  fuzzy $\mathscr{F}$-metric spaces.  Section 4 contains a type new fixed-point theorems under generalized contraction conditions involving $\psi$-functions and show the  applicability of the obtained  result  by solving the satellite web coupling
problem.

\section{Preliminaries}
      In order to carry this study, we need to provide the following definitions and results. 
      
      \begin{dfn} \cite{1} \label{dfn 6}
      	A binary operation $ \star : [0 , 1] \times [0 , 1] \rightarrow [0 , 1] $ is called a  $t$-norm if it satisfies the following conditions:
      	\begin{enumerate}[(i)]
      		\item  $ * $ is associative and commutative; 
      		\item $ \alpha ~ * ~ 1 = \alpha, ~ \forall \alpha \in [0 , 1] $;  
      		\item $ \alpha * \gamma \leq  \beta * \delta ~ $  whenever $ \alpha \leq  \beta  $ and $  \gamma \leq   \delta,  ~\forall \alpha, \beta,
      		\gamma, \delta \in [0 , 1] $.
      	\end{enumerate}
      	If $ * $ is continuous then it is called continuous $t$-norm.
      \end{dfn}

      \begin{eg} \cite{1} 
      	The following are examples of some t-norms given in \cite{1}. 
      	\begin{enumerate}[(i)]
      		\item  Standard intersection: $ \alpha *   \beta = \min \{ \alpha ,   \beta \} $.
      		\item   Algebraic product: $ \alpha *   \beta = \alpha    \beta $.
      		\item Bounded difference: $ \alpha *   \beta = \max \{ 0,  \alpha +  \beta -1 \} $.
      	\end{enumerate}
      \end{eg}

   Following is the definition of    fuzzy metric space given by George \& Veeramani \cite{G3}.

      \begin{dfn} \cite{G3} \label{dfn 7}
      	The $ 3$-tuple $ (X, M, \star ) $ is said to be a (George \& Veeramani type) fuzzy metric space if  the  fuzzy set $ M $   on $ X^2 \times (0, \infty ) $ satisfies the following conditions:
      	\begin{enumerate}[($ M 1 $)]
      		\item  $ M (x,y,t) > 0 $,
      		\item $ M (x,y,t) = 1 $ if and only if $ x = y $, 
      		\item $ M (x,y,t) = M (y, x, t) $
      		\item $ M (x,y,t) \star M (y, z, s) \leq M (x, z, s+t) $,
      		\item  $ M (x, y, .) : (0, \infty ) \to ( 0, 1] $ is continuous, 
      	\end{enumerate}
      	for all $ x, y, z \in X $ and $ t, s > 0 $.
      \end{dfn}

      \begin{lma} \cite{2}  \label{ lma 1}
      	George \& Veeramani  fuzzy metric $ M (x, y, \cdot ) $ over a nonempty set $X$ is non-decreasing with respect to $ t > 0 $ for all $ x, y \in  X $. 	
      \end{lma}

      \begin{rem}
      	\label{important rem}
      		In  fuzzy metric space $ (X, M, \star ) $, George \& Veeramani considered the followings.
      	\begin{enumerate}[(i)] 
      		\item If $ M (x, y, t) > 1 - r $ for all $ x, y \in  X $, $ t > 0 $, $ 0 < r < 1 $, we can find a $ t_0, 0 < t_0 < t $ such that $ M (x, y, t_0) > 1 - r $. 
      		\item For any $ r_1 > r_2 $ in $ ( 0, 1 )$, we can find  $ r_3 \in (0,1) $ such that $ r_1 \star r_3 \geq r_2 $ and for any $ r_4 \in (0,1) $ we can find  $ r_5 \in (0,1) $ such that $ r_5 \star r_5 \geq r_4  $.
      	\end{enumerate}
      \end{rem}

 Next we recall the definition of  fuzzy $ \mathscr{F} $- metric space introduced by Das et al.  \cite{4}. 
  Throughout the article,	$\mathscr{F}$ denotes  the set of all functions $ \textit{f}: [0,1] \to [0,1] $ which satisfy the following conditions: 
  \begin{enumerate}[($\mathscr{F}1$)]
  	\item \textit{f} is strictly increasing in $ [0,1) $;
  	\item For every sequence $\{ t_n\} $ in $ [0,1] $ we have,
  	$\underset {n \to \infty}{\lim} ~ t_n = 1 \iff \underset {n \to \infty}{\lim}  ~ \textit{f}(t_n)=1. $ 
  \end{enumerate}

  \begin{eg} \cite{4} \label{ eg 1}
  	Following are some examples of members of $\mathscr{F}$. 
  	\begin{enumerate}[(a)]
  	\item  $ \textit{f} ( x ) = x^n, ~ \forall x \in [ 0, 1 ] $, $ n\in \mathbb{N}$.  
    \item $  \textit{f} ( x ) = \sqrt{x}, ~ \forall x \in [ 0, 1 ] $.
  	\end{enumerate}
  \end{eg}

  \begin{dfn}\cite{4} \label{main dfn}
  	Let $ X $ be a non empty set and $ M : X \times X \times  (0, \infty ) \to [0,1] $ be a mapping and $\star$ be a continuous t-norm. 
  	If there exists $ ( \textit{f}, \alpha) \in \mathscr{F} \times (0,1] $ such that $ M $ satisfies the following conditions: 
  	\begin{enumerate}[($\mathscr{F}M1$)]
  		\item $ M(x,y,t) > 0  $; $ ~~ \forall x, y \in X $ \& $ t > 0 $;
  		\item $ (  M(x,y,t) = 1, ~ \forall t > 0  ) $  iff $ x = y $;
  		\item $ M(x,y,t) = M(y,x,t)  $; $ ~~ \forall x, y \in X $ \& $ t > 0 $;
  		\item for every $ (x, y) \in X \times X  $, for every $ N \in \mathbb{N} $, $ N \geq 2$  and for every $ \{ u_i\}_i ^N \subseteq  X $ with $ u_1 = x $ and $ u_N =  y $, we have 
  		$$ M(x,y,t) < 1 \implies ( \textit{f}(M(x,y,t)))^\alpha \geq \textit{f}( M ( u_1, u_2 ,t_1 ) ~ \star ~ M(u_2, u_3, t_2) ~ \star  \cdots  \star ~ M(u_{N-1}, u_N, t_{N-1} ) ) $$
  		where $t = t_1 + t_2 + \ldots + t_{N-1} $; $ t_i > 0 $ for $  i = 1, 2, \cdots, (N - 1) $,  
  	\end{enumerate}
  	then $  M $ is said to be a fuzzy $\mathscr{F}$- metric on $ X $ and the 5-tuple $ ( X, M, \mathscr{F}, \alpha,  \star) $ is said to be a fuzzy $ \mathscr{F} $-metric space.
  \end{dfn}

  \begin{eg} \cite{4} \label{ex from 1st paper}
  	Let $ X = \mathbb{R} $ and define a function $ M : X \times X \times (0,\infty) \rightarrow [0,1]  $ by \\
  	$$ M (x,y,t) = \left( \dfrac{t}{t+1}\right)^{|x-y|^2} ~~  \text{for all}~  (x,y) \in X \times X ~ \text{and} ~ t \in (0,\infty). $$ 
  	Then $ M $ is a fuzzy $ \mathscr{F} $-metric on $X$ with respect to the `product' t-norm.
  \end{eg}

  \begin{ppn} \cite{4} \label{conv and Cauchy}
  	Let $(X, M, \textit{f}, \alpha, \star )$ be a fuzzy $\mathscr{F}$-metric space, $\{x_n\} \subseteq X $ be a sequence and $ x \in X$. Then
  	\begin{enumerate} [(i)]
  		\item $\{ x_n \} $ is $\mathscr{F}$-convergent to $ x $ iff    {$\underset{n \to \infty} {\lim}  M(x_n, x, t) = 1 ~~ \forall~ t > 0 $.}
  		\item $ \{x_n \} $ is $\mathscr{F}$-Cauchy iff   $\underset{m, n \to \infty}{\lim} M(x_n, x_m, t) = 1 ~~ \forall~ t > 0$.
  	\end{enumerate}
  \end{ppn}

\begin{ppn} \cite{4}
	Limit of an $ \mathscr{F} $-convergent sequence in $ ( X, M, \textit{f}, \alpha,  \star) $ is unique.
\end{ppn}

\begin{ppn}\label{proposition 2.11 ( 9.01.2025)}
    \cite{4} In a fuzzy $ \mathscr{F}$-metric space, every $ \mathscr{F}$-convergent sequence is an $ \mathscr{F}$-Cauchy sequence.
\end{ppn}
 \begin{ppn}\label{ppn 7}
 	In a fuzzy $\mathcal{F}$-metric space $(X, M, \textit{f}, \alpha,\star)$, for any subset $A \subset X $, 
 	$x \in$ ${\bar{A}}$, $ 1>r>0 $  implies   $B_\mathcal{F}(x,r,t)$ $\cap A \ne$ $\phi$ for any $ t> 0$.
 \end{ppn} 	
\begin{proof}
 		Suppose $x \in$ ${\bar{A}}$ and $ 1>r>0 $ be  fixed.\\Then $\exists$ a sequence $\{x_n\}$ $\subseteq A$ such that $\underset{n \to \infty }{\lim}$ $M(x_n, x,t) = 1$ for any $ t>0 $
 		\begin{align*}
 			&\implies \exists N \in \mathbb{N} ~ such ~ that ~ M(x_n, x,t) > 1-r  ~\forall~ n \geq N\\
 			&\implies x_n \in B_\mathcal{F}(x,r,t)~ \forall ~n \geq N\\
 			&\implies B_\mathcal{F}(x,r,t) \cap A \ne \phi .
 		\end{align*}
 	\end{proof}

We develop a fixed point theorem in Section 4 using the following class of mappings. 

\begin{dfn} \cite{3} \label{dfn of si function}
	Let $ \Psi $ be the class of all mappings $ \psi : [ 0, 1] \to [0, 1]  $ satisfying
\begin{enumerate}[(i)]
	\item $ \psi $ is continuous and nondecreasing;
	\item $ \psi ( t ) > t ~~ \forall t \in ( 0, 1 ) $.
\end{enumerate}	 
\end{dfn}

\begin{eg}\cite{3} 
Define a function $ \psi  : [0, 1] \to [0, 1] $ by $ \psi ( t ) = \dfrac{t}{t + k ( 1 - t ) }, ~ t \in [ 0, 1 ] $ where $ k \in ( 0, 1 ) $. Then $ \psi  \in \Psi $. 
\end{eg}

\begin{lma}\cite{3} \label{prop of si 1}
	 If $ \psi \in \Psi $ then  $ \psi ( 1 ) = 1 $. 
\end{lma}

\begin{lma}
	\cite{3} \label{prop of si 2}
	If $ \psi\in \Psi $ then    $ \underset{ n \to \infty}{\lim} \psi^n ( t ) = 1 $  for all $ t \in ( 0, 1 ) $. 
\end{lma}

\begin{thm}\cite{3}
	Let $ (X, M, \star) $  be a fuzzy metric space	satisfying the $ \psi $-contraction condition 
	$$ M(x, y, t) > 0  \iff M ( S (x), S (y), t) \geq \psi (M(x, y, t)) ~~ \forall t > 0 $$ 
	where $ S :	X \to X $ is   a   self mapping in $X$. Then $ S $	has at most one fixed point.
\end{thm}

\section{Some basic properties of fuzzy $ \mathscr{F} $-metric space}
        
   In this section, we develop   fundamental results concerning the notions of completeness, compactness and total boundedness, which are essential to study  the structural properties of the concerned fuzzy $ \mathscr{ F } $-metric space.

    \begin{lma}
  Suppose  $ f : [ 0, 1 ] \to [ 0,1 ] $ be a function satisfying the following conditions: 
     \begin{enumerate}[(i)]
       \item $ f $ is strictly increasing in $ [ 0, 1 ) $.
          \item For every sequence $ \{ t_n \} $ in $ [ 0, 1 ] $ we have,  $ \underset{ n \to \infty } {\lim} t_n = 1 \iff \underset{ n \to \infty }{\lim} f ( t_n ) = 1 $.   
         \end{enumerate}
                  Then $ f ( 1 ) = 1 $.
    \begin{proof}
                  	Suppose that $ f ( 1 ) \ne 1 $.
                  	So, $ f ( 1 ) < 1 $ or $ f ( 1 ) > 1 $.
                  	
                  	\textbf{ Case I : } 
                  	 If $ f ( 1 ) > 1 $, then $ f ( 1 ) \notin [ 0, 1 ] $. 
                  	So, this case is not possible.
                  	
                  	\textbf{ Case II : }
                  	 If $ f ( 1 ) < 1 $.\\
                  	 Let us consider a sequence $ \{ t_n \} $ by $ t_n = 1 $ for all $ n \in \mathbb{ N } $. 
                  	 So, by condition $ (ii) $,   $ \underset{ n \to \infty }{\lim}  f ( t_n ) = 1 $. 
                  	 Again, 
                  	 \begin{align*}
                  	 	& t_n = 1  ~ ~~ \forall ~  n   \in \mathbb{ N }\\
                  	 	 \implies & f ( t_n ) = f( 1 ) ~ ~~ \forall ~ n \in \mathbb{ N }\\
                  	 	\implies & \underset{ n \to \infty } {\lim} f ( t_n ) = f ( 1 ) < 1, ~~~\forall ~ n \in \mathbb{ N } \\
                  	 	\implies & \underset{ n \to \infty } {\lim} f ( t_n ) < 1, \text{contradiction to the condition (ii)}.  
                  	 \end{align*}  
                   So, our assumption is wrong. 
                   Hence, $ f ( 1 ) = 1 $.
                  \end{proof}
                  \end{lma}

        \begin{ppn}
            If $ \{  x_n \}$ be a $ \mathscr{F} $-convergent sequence in $ ( X, M, f, \alpha, \star ) $ converging to $ x $. Then $ \{  x_{n_k} \} $ be a subsequence which is also $ \mathscr{F} $-convergent converging to the same point $ x $.
        \end{ppn}
            \begin{proof}
                Since $ \{  x_n \} $ is a $ \mathscr{ F } $-convergent sequence converging to $ x $, then $$ \underset{ n \to  \infty } {\lim}  M (  x_n, x, t ) = 1 ~ \forall t > 0.  $$
                
Let $ \{  x_{ n_k} \}$ be a subsequence of $ \{ x_n \}$. 
By  proposition \ref{proposition 2.11 ( 9.01.2025)}, $ \{  x_n \} $ is a $ \mathscr{ F} $- Cauchy sequence.  
Therefore,  
$$ \underset{ n \to  \infty }  {\lim}  M (  x_{n_k}, x_n, t ) = 1 ~ \forall t > 0 $$.

Now, 
$$ \underset{ n \to \infty }  {\lim}   \{  M (  x_n, x, t ) \star M ( x_n, x_{n_k}, t ) \} = 1. $$
 
Again now, we have 
\begin{align*}
    & \left(  f ( M ( x_{n_k}, x, 2t)\right)^\alpha   \geq f ( M ( x_{n_k}, x_n, t) \star M ( x_n, x, t ) )\\
    & \underset{ n \to \infty } {\lim }  \left(  f ( M ( x_{n_k}, x, 2t)\right)^\alpha \geq   \underset{ n \to \infty } {\lim}     f ( M ( x_{n_k}, x_n, t) \star M ( x_n, x, t ) ) =  \underset{ n \to \infty } {\lim}   f ( 1 \star 1 ) = 1~ ( \text{ using } \mathscr{ F}2).  
\end{align*} 

Therefore,
\begin{align*}
    &  \left(  f ( M ( x_{n_k}, x, 2t)\right)^\alpha  = 1  ~ \forall t> 0\\
   & \implies M (  x_{n_k}, x, 2t ) = 1  ~ \forall t> 0.
    %%%%%%%%%%%%
\end{align*}
This shows that $ \{ x_{n_k} \} $ converges to $ x $. 
Hence, the proof is complete.
%%%%%%%%%%%%
\end{proof}

          %%%%%% every compact subset of a fuzzy F-matric sp ace is closed and bounded.
        \begin{dfn}
           	Let $ ( X, M,  f, \alpha, \star ) $ be a fuzzy $ \mathscr{F} $-metric space and $ A \subset X $. $ A $ is said to be compact if every sequence in $ A $ has a convergent subsequence which converges to some point in $ A $.
        \end{dfn}

                \begin{thm}
           	  Every compact subset of a fuzzy $ \mathscr{F} $-metric space is closed and bounded.
              \end{thm}
           	 \begin{proof}
           	 	 %%%
           	 	 Let $ ( X, M,  f, \alpha, \star ) $ be a fuzzy $ \mathscr{F} $-metric space and $ A $ be a compact subset of $ X $.  
           	 	 If possible suppose that $ A $ is not closed. So there exists a sequence $ \{ x_n\} $ in $ A $ such that $ x_n \to x $ as $ n \to \infty $ but $ x \notin A $. 
                 
           	 	 Since $ A $ is compact, so there exists a subsequence $ \{ x_{n_k}\} $ of $ \{ x_n\} $ which converges to some point in $ A $. 
                 
           	 	 But $ x_n \to x $ as $ n \to \infty $ implies   $ x_{n_k} \to x $ as $ n \to \infty $  and hence $ x \in A $.  
           	 	 Which is a contradiction.  
           	 	 Thus $ A $ is closed. 
                 
           	 	 %%%%
           	 	 Now we show that $ A $ is bounded. If possible suppose  that $ A $ is unbounded. 
           	 	 Let $ \{ x_n\} $ be a sequence in $ A $.   So it has a convergent subsequence $ \{ x_{n_k}\} $ which converges to some $ x \in A $.
               Let $ x_0 \in A $ be a fixed element. 
               
           	 	 Now choose $ 0 < \epsilon < 1$. Then by $ ( \mathscr{F}2 ) $,  there exists $ \delta \in ( 0,1) $ such that 
           	 	 %%%%
           	 	 \begin{equation}\label{eq 1, th 3}
           	 	 	1-\delta < t \leq 1 \implies 1-\epsilon < f(t) \leq 1.
           	 	 \end{equation} 
            	 %%%%%%
                 
            	 Again there exists $ t_0 > 0 $ such that $ M ( x_0, x, t_0 ) > ( 1 - \delta ) $. 
            	 By Remark \ref{important rem}, we can find $ r \in (0,1) $ such that
            	 %%%%%
            	 \begin{equation}\label{eq 2, tm 3}
            	 	M ( x_0, x, t_0 ) \star (1-r) > (1-\delta)
            	 \end{equation}
             %%%%%
             
               We choose $ \{ \alpha_n \} \in (0,1) $ such that $ \alpha_n \to 1 $ as $ n \to \infty $.  
               So for   given $ t'  (> t_0 )  $, for each $ \alpha_k $ there exists $ x_{n_k} $ in $ A $ such that 
               %%%%%%
               \begin{equation}\label{eq 3, th 3}
               	 M ( x_{n_k}, x_0, t'   ) \leq ( 1 - \alpha_k ).
               \end{equation}
           %%%%%%%%%%
           
              Since $ x_{n_k} \to x $ as $ n \to \infty $, so for $ ( t' - t_0) = t_1 $(say), there exists $ m (t_1) \in \mathbb{N} $ such that 
             $$ M ( x, x_{n_k},  {t_1} )  > (1-r) ~~ \forall ~ k \geq m (t_1 ). $$ 
             
              Hence we have,
              %%%%%%%%%%%%
              $~ M ( x_0, x, t_0) \star M ( x, x_{n_k}, t_1) >  M ( x_0, x, t_0) \star (1-r) > (1- \delta) ~~~ \forall ~ k \geq m (t_1) ~~~  ( using ~ (\ref{eq 2, tm 3})) ~ $ 
             %%%%%
             which implies
             %%%% 
             \begin{align*} 
                &   f ( M ( x_0, x, t_0) \star M ( x, x_{n_k}, t_1) ) > 1- \epsilon ~~~~~~~~~~~ \forall ~ k \geq m (t_1) ~~~~~~~~~ ( using ~ ( \ref{eq 1, th 3})) \\
                  \implies & \left( f ( M ( x_0, x_{n_k}, t') ) \right)^\alpha > 1-\epsilon ~~~~~~~~~~~ \forall ~ k \geq m (t') ~~~~~~~~~ ( using ~ ( \mathscr{F}M4 ) ) \\
             	  \implies & \left( f( 1-\alpha_k) \right)^\alpha > 1-\epsilon  ~~~~~~~~~~~ \forall ~ k \geq m  ~~~~~~~~~ ( using ~ (\ref{eq 3, th 3})) \\
             	 \implies & \underset{k \to \infty }{\lim } \left( f( 1-\alpha_k) \right)^\alpha  \geq 1-\epsilon. 
             	\end{align*}
               Since $ 0 < \epsilon < 1 $  is  arbitrary, we can write
             	 \begin{align*}
             	    \underset{k \to \infty }{ \lim } \left( f( 1-\alpha_k) \right)^\alpha = 1  
             	   ~ or ~ \underset{k \to \infty }{ \lim }  f( 1-\alpha_k) = 1  
             	   ~ or ~ \underset{k \to \infty }{ \lim } ( 1-\alpha_k) = 1 ~  ~( \text{using} ~ (\mathscr{F}2)) 
             	   ~ or ~ \underset{k \to \infty }{ \lim } \alpha_k = 0.
             \end{align*}  
              This is a contradiction to our assumption. Hence $ A $ is bounded.
           	 \end{proof}

           \begin{thm}
           	In a fuzzy $ \mathscr{F} $-metric  space, every compact set is complete.
          \end{thm}
           	\begin{proof}
           		Let $ A \subseteq X $ be a non-empty compact  set in a fuzzy $ \mathscr{F} $-metric  space $ ( X, M,  f, \alpha, \star ) $ . 
           		Let $ (f,\alpha) \in \mathscr{F} \times (0,1] $ be such that $ (\mathscr{F}M4) $ holds in $ (X, M,  f, \alpha, \star ) $.

           		Let $ 0 < \epsilon < 1 $. Then by $ (\mathscr{F}2) $ there exists $ 0 < \delta < 1 $ such that 
           		\begin{equation}\label{thm 5, eq 1}
           			1-\delta < t \leq 1 \implies 1-\epsilon < f(t) \leq 1.
           		\end{equation} 
           	For $ \delta \in ( 0,1 ) $, there exists $ \lambda \in (0,1) $ such that $ ( 1 - \lambda )  > ( 1 - \delta ) $. Then using
           	   Remark $ \ref{important rem} $, for   $ \lambda \in ( 0,1 ) $,  we can find $ r_1 \in (0,1) $ such that $ ~
           	   	( 1-r_1) \star ( 1-r_1) \geq (1- \lambda ) > ( 1 - \delta ). $

              Again using Remark $ \ref{important rem} $, we can find $ r_3 \in (0,1) $ 
              \begin{equation*}
              	( 1-r_1) \star ( 1-r_1) \star (1- r_3) \geq (1 - \delta). 
              \end{equation*}
          %%%%%%%
            Let $ \{ x_n\} $ be an $ \mathscr{F} $- Cauchy sequence in $ A $.  
            Then for $ 0 < r_1 < 1 $ and for each $ t_1 > 0 $, there exists $  N_1 ( t_1) \in \mathbb{N} $ such that  
             $ M ( x_n, x_m, t_1) > 1-r_1, \quad \forall m,n \geq  N_1 ( t_1) $.  
             
            Since $ A $ is compact, so  $ \{ x_n\} $ has a $ \mathscr{F} $-convergent subsequence, say $ \{ x_{n_k} \} $ which converges to $ x $ in $ A $.   
            Then for $ 0 < r_3 < 1 $ and for each $ t_2 > 0 $ there exists $ N_2(t_2)  \in \mathbb{N} $ such that 
            $ M (x_{n_k}, x, t_2) > 1-r_2, ~ \forall k \geq N_2 (t_2) $. 
            
            Let $ T = 2 t_1  + t_2 $ and $ N(t) =  \text{max}  \{ N_1(t_1), N_2(t_2)   \} $.   
            Then,
            \begin{align*}
            &	M (x_n, x_m, t_1) \star M ( x_m, x_{n_k}, t_1) \star M ( x_{n_k}, x, t_2) > (1-r_1) \star (1-r_1) \star (1-r_3) \quad \forall~ m,n,k \geq N(t) \\
             \implies &	M (x_n, x_m, t_1) \star M ( x_m, x_{n_k}, t_1) \star M ( x_{n_k}, x, t_2) > 1-\delta \quad \forall~ m,n,k \geq N(t) \\
              \implies & f \left( 	M (x_n, x_m, t_1) \star M ( x_m, x_{n_k}, t_1) \star M ( x_{n_k}, x, t_2) \right) > 1-\epsilon \quad \forall~ m,n,k \geq N(t) \quad ( \text{using}~ (\ref{thm 5, eq 1}) )\\
             \implies & \left( f ( M (x_n, x, T ) ) \right)^\alpha > 1-\epsilon \quad \forall~ m,n,k \geq N(t).  
            \end{align*}

            Since $ 0 < \epsilon < 1 $ is arbitrary, so letting $ n \to \infty $, we get 
            \begin{equation*}
            	  \underset{n \to \infty}{ \lim } f \left( M ( x_n, x, T )\right) = 1 \quad  {\forall~ T > 0}  \\
            	\implies \underset{n \to \infty }{ \lim} M ( x_n, x, T ) = 1 \quad ( \text{using} ~ ( \mathscr{F}2 ) ) 
            \end{equation*}
            
         This proves that  $ \{ x_n\} $ is a $ \mathscr{F} $-convergent sequence in $ A $ converging to $ x \in A $ thus $ A $ is complete in $ X $. 
           	\end{proof}

                %%%%%%%%%%%%%%  r- totally bounded
         \begin{ppn}\label{ppn 1}
          Every finite subset in a fuzzy $ \mathscr{F} $-metric space is bounded.
          \end{ppn}
          \begin{proof}
          	Let $ ( X, M,  f, \alpha, \star ) $ be a fuzzy $ \mathscr{F} $-metric space and $ A $ be a finite subset of $ X $ containing $ n $ elements $ x_1, x_2, \ldots x_n $.  
          	Choose $ t_0 > 0 $ fixed.  
          	Let $ \underset{i, j}{ min } ~ M (x_i, x_j, t_0 ) = \beta $.  
          	Clearly $ \beta \in (0,1) $.  
          	Then we can choose $ r \in (0,1) $ such that $ \beta > 1-r $  
          	$ \implies M (x_i, x_j, t_0 ) > 1-r \quad \forall~ x_i, x_j \in A $  
          	$ \implies  A $ is bounded.
          \end{proof}

        \begin{dfn}
        	Let $ ( X, M,  f, \alpha, \star )$ be a fuzzy $ \mathscr{F} $-metric space and $ A \subset X $ and $ r \in (0,1) $ be given. Let $ \epsilon > 0 $ be a positive number. A subset $ B \subset X $ is said to be an $ r-\epsilon $-net for the set $ A $ if for any $ x \in A$  there exists $ y \in B $ such that $ M ( x, y, \epsilon ) > 1-r $.  
        \end{dfn}

    % $ B $ may be finite or infinite.

     \begin{dfn}
     	A set $ A $ in a fuzzy $ \mathscr{F}$-metric space $ ( X, M,  f, \alpha, \star) $ is said to be $ r $-totally bounded  if for any $ \epsilon > 0 $ there exists a finite $ r- \epsilon $-net for the set $ A $.
     \end{dfn}

        \begin{thm} \label{totally bdd implies bdd}
        		Let $ ( X, M, f, \alpha, \star ) $ be a fuzzy $ \mathscr{F} $-metric space and $ A \subset X $. If $ A $ is  $ r $-totally bounded for some $ r \in (0, 1) $ then $ A $ is bounded.
         \end{thm}		
        	\begin{proof}
        	Since $ A $ is  $ r $-totally bounded, so for any $ \epsilon > 0 $, there exists a finite $ r-\epsilon $-net $ B $ for $ A $.
         Since $ B $ is finite, thus $ B $ is bounded. So there exists $ t_0 > 0 $ and $ r_0 \in ( 0, 1 ) $ such that
        	%%%%%%%%%%%%%
        	\begin{equation}\label{Eq 2, thm 1}
        		M ( y_1, y_2, t_0 ) > 1-r_0 \quad \forall ~ y_1, y_2 \in B
        	\end{equation} 
        	%%%%%%%%%%%%%%%
Again  for each $ a \in A $,  there exists $ b \in B  $ such that 
        	%%%%%%%%%%%
        	\begin{equation}\label{Eq 1, thm 1}
        		M ( a, b, t_0 ) > 1-r	
        	\end{equation}
        	%%%%%%%%%%%%%
        	Now choose $ a_1, a_2 \in A $ arbitrarily. Then there exist  $ b_1, b_2 \in B $ such that  
        	%%%%%%%%%%%%%%%%
        	\begin{equation}\label{Eq 3, thm 1}
        		M ( a_1, b_1, t_0 ) > 1-r  \quad \text{\&} \quad
        		M ( a_2, b_2, t_0 ) > 1-r
        	\end{equation}
        	%%%%%%%%%%%%%%%
        	Again using $ ( \mathscr{F}M4 ) $, we have 
        	%%%%%%%%%%%%%%
        	\begin{equation}\label{Eq 4, thm 1}
        		\left( f ( M ( a_1, a_2, t ) ) \right)^\alpha \geq f ( M (a_1, b_1, t_0 )  \star M ( b_1, b_2, t_0) \star M ( b_2, a_2, t_0 ) ) \quad \text{ where } t = 3 t_0 ~(\text{fixed}).
        	\end{equation}
        	%%%%%%%%%%%%
        	Now from  $ (\ref{Eq 2, thm 1}) ~ \& ~ (\ref{Eq 3, thm 1}) $ we have 
        	%%%%%%%%%%%%%%%
        	\begin{equation}\label{eq 5, thm 1}
        		M (a_1, b_1, t_0 )  \star M ( b_1, b_2, t_0 ) \star M ( b_2, a_2, t_0 ) > ( 1-r ) \star ( 1-r_0 ) \star ( 1-r )
        	\end{equation}
        	%%%%%%% 
        	If $	M (a_1, b_1, t_0 )  \star M ( b_1, b_2, t_0 ) \star M ( b_2, a_2, t_0 ) = 1 $ then
        	%%%%%%%%%%%%%%
        	\begin{align*}
        		&\left( f ( M ( a_1, a_2, t ) ) \right)^\alpha = 1 \quad \forall ~ a_1, a_2 \in A  \\
        			\implies & M ( a_1, a_2, t ) = 1  \quad \forall ~ a_1, a_2 \in A \quad ( \text{using} ~ ( \mathscr{F}2 ) ) \\
        		 \implies & M ( a_1, a_2, t ) > (1-\beta) \quad \forall ~ ~ a_1, a_2 \in A ~ \text{and for any } \beta \in ( 0, 1 )
        	\end{align*}
        	This shows that $ A $ is bounded. \\
        	%%%%%%%%%%%%%
        	If $ M (a_1, b_1, t_0 )  \star M ( b_1, b_2, t_1 ) \star M ( b_2, a_2, t_0 ) < 1 $  
        	then 
        	\begin{align*}
        	&	f \left( M (a_1, b_1, t_0 )  \star M ( b_1, b_2, t_1 ) \star M ( b_2, a_2, t_0 ) \right) < 1  \\
        		or ~~ &  f \left( (1-r) \star (1-r_0) \star (1-r) \right) < 1 \quad  ( using  (\mathscr{F}2 ))
        	\end{align*}
        	So we may write, 
        $ f \left( (1-r) \star (1-r_0) \star (1-r) \right) = (1-\gamma), \quad \text{for some} ~ \gamma \in (0,1). $ \\
        	Thus from $ (\ref{Eq 4, thm 1}) $ we have, 
        	\begin{align}\label{modification 1}
        	   & \left( f ( M ( a_1, a_2, t ) ) \right)^\alpha > (1-\gamma), \quad \forall ~ a_1, a_2 \in A \nonumber \\ 
             or ~~ &  f ( M ( a_1, a_2, t ) ) > (1-\gamma )^\frac{1}{\alpha } = (1-\delta) (\text{say}) \quad \forall~ a_1, a_2 \in A , \text{ where } \delta \in (0,1).
        	\end{align}  
          %%%%%%%%%%%%%%
Using $ ( \mathscr{F}2)$, for $ \delta \in ( 0, 1 )  $ there exists $ \mu \in ( 0, 1 ) $ such that 
$$ 1 - \mu < t \leq 1  \iff 1 -\delta < f ( t) \leq 1. $$
Therefore the relation (\ref{modification 1}) gives  $ ~ M (a_1, a_2, t ) > (1- \mu ) ~ ~ \forall ~ a_1, a_2 \in A $.  
 Since $ t $ is fixed, so it follows that $ A $ is bounded.	
        	\end{proof}

    \begin{rem}
    	The converse of the above Theorem \ref{totally bdd implies bdd} is not necessarily true. For justification we consider the following example.
    \end{rem}         
    
    \begin{eg} 
        Recall the fuzzy $ \mathscr{F} $-metric space of Example \ref{ex from 1st paper}.  	 Consider the subset    $ A = \{ 1, 3, 5,7, 9 \} $  of $ \mathbb{ R } $. 
            	Then by Lemma \ref{ppn 1}, $ A $ is bounded set.\\
            	If possible suppose that $ A $ is $ r $-totally bounded.\\
            	So, for $ \epsilon = \frac{ \sqrt{ 2 }}{ 3 }, r = 1- \frac{ 1 }{ \sqrt{ 2 } } $ there exists a finite $ r - \epsilon $-net, say $ N $ for $ A $. Thus
            	for $ x_i, x_j ( i \ne j ) \in A $, there exist  $ y_i, y_j \in N $ such that
            	  $$ M ( x_i, y_i, \epsilon ) > (1 - r) \text{ and } M ( x_j, y_j, \epsilon ) > (1-r).   $$
            	 Now,
            	 \begin{align*}
            	 ( f ( M ( x_i, x_j, 2 \epsilon ) ) )^\alpha & \geq f ( M ( x_i, y_j, \epsilon ) \star M ( y_j, x_j, \epsilon ) ) \\
            	 & \geq f ( ( 1-r) \star ( 1-r ) ) ~~~~ (\text{since} ~ f ~ \text{is an increasing function})\\
            	 & = f ( ( 1-r )^2 )  
            	  \geq f ( \frac{ 1 }{ 2 } )  = \frac{ 1 }{ \sqrt{ 2 } }.
            	\end{align*}
              Again  
              \begin{align*}
              	\left( M ( x_i, x_j, 2 \epsilon ) \right)^{ \alpha } & =  \left( \frac{ 2 \epsilon }{ 1 + 2 \epsilon } \right)^{ \frac{ | x_i - x_j |^2  }{ 2 } } 
              	 = \left( \frac{ 2 \frac{ \sqrt{ 2 } }{ 3 } }{ 1 + 2 \frac{ \sqrt{ 2 }}{ 3 } } \right)^{ \frac{ | x_i - x_j |^2 }{ 2 } }  
              	  = ( 0.49 )^{ \frac{ | x_i - x_j |^2 }{ 2 } }
              \end{align*}
           If $ x \ne y $ then $ \frac{ | x_i - x_j |^2 }{ 2 } \geq 2 $. 
           Therefore, $ ( 0.49 )^{ \frac{ | x_i - x_j |^2 }{ 2 } } < \frac{ 1 }{\sqrt{ 2 } }  $. It is a contradiction.\\
           So, there is no $ r- \epsilon $-net for $ A $.
           Hence $ A $ is not $ r $-totally bounded.	 
    \end{eg}

         \begin{thm}
         	Let $ ( X, M,  f, \alpha, \star ) $ be a fuzzy $ \mathscr{F} $-metric space and $ A \subset X $. If $ A $ is a compact set then $ A $ is $ r $-totally bounded for all $ r \in (0,1) $.
          \end{thm}
         	%%%%%%%%%%%%
         	\begin{proof}
         		We assume that $ A $ is compact. We
         		choose fixed $ r \in (0,1) $ and $ \epsilon > 0 $   arbitrarily. Then by $ ( \mathscr{F}2 ) $ there exists $ 0 <  \delta < 1 $ such that
         		%%%%%%%%%%%%%
         		\begin{equation}\label{th 10,eq 1}
         			1-\delta < t  \leq 1 \implies 1-\epsilon < f(t) \leq 1.
         		\end{equation}
         		For $ \delta \in (0,1) $, we can find $ r_1 \in (0,1) $ such that 
         		%%%%%%%%%%%%%%
         		\begin{equation}\label{thm 10, eq 2}
         			(1 - r_1 ) \star ( 1 - r_1 ) \geq ( 1 - \delta ). 
         		\end{equation}
         		Let $ x_1 $ be an arbitrary element of $ X $. If $ M ( a, x_1, \epsilon ) > ( 1-r) ~ \forall ~ a \in A $,  
         		then a finite $ r-\epsilon $-net $ B $ exist for $ A $ i.e. $ B = \{x_1\} $.
         		If not there exists a point $ x_2 \in A $ such that $ M ( x_1, x_2, \epsilon ) \leq ( 1 - r ) $. \\
         		If  for all $ a\in A $, $ M ( a, x_1, \epsilon ) > ( 1-r) $ or $ M ( a, x_2, \epsilon ) > ( 1-r) $ then $ B = \{ x_1, x_2 \} $ is a finite $ r-\epsilon $-net  exists for $ A $. 
         		Continuing in this way, we obtain points  $ \{x_1, x_2, \ldots, x_n \}$ where $ x_1 \in X $, $ x_2, x_3, \cdots x_n \in A $ such that   $ M (x_i, x_j, \epsilon ) \leq ( 1 - r ) $ for $ i \ne j $.  Here two  cases may arise.

         		 \textbf{Case-I} The procedure stops after $ k^{th } $- steps. 
                 
         		 Then we obtain points $ x_1, x_2, x_3, \ldots, x_k $ such that for every $ a \in A $ at least one of the inequalities $ M (x_i, a, \epsilon) > (1-r), i= 1,2, \ldots, k $ holds and then $ B = \{ x_1, x_2, \ldots, x_k \} $ is a finite $ r-\epsilon $-net for $ A $. 
         		 Hence $ A $ is $ r $-totally bounded. 

              \textbf{Case-II}  
              The procedure continues infinitely.  
              Then we obtain an infinite sequence $ \{ x_n \} $ where $ x_1 \in X  $ and $ x_i \in A $ for $ i \geq 2 $ such that 
              \begin{equation}\label{thm 10, eq 2,5}
              	 M (x_i, x_j, \epsilon ) \leq ( 1-r)  ~~~ \text{for} ~  i \ne j.  
              \end{equation} 
              Since $ A $ is compact then there exists a subsequence $ \{ x_{n_k} \} $ of $\{ x_n \} $ which converges to $ x \in A $.
               Then for   $ \frac{\epsilon}{2} > 0 $ and $ r_1 \in (0,1) $ there exists $ N (\epsilon) \in \mathbb{N} $ such that  
               \begin{equation}\label{thm 10, eq 3}
               	M ( x_{n_i}, x, \frac{\epsilon}{2}   )  > (1-r_1) ~  ~~ \forall ~ i \geq N.
               \end{equation}

               Hence we have,
             \begin{equation}\label{thm 10, eq 3.5}
             	 M (x_{n_i}, x, \frac{\epsilon}{2} ) \star  M (  x_{n_{i+1}}, x,  \frac{\epsilon}{2} ) > ( 1 - r_1 ) \star ( 1 - r_1 ) \geq ( 1 - \delta ) ~  ~~ \forall ~ i \geq N.
             \end{equation}
             %%%
             
               Again using $ (\mathscr{F}4) $, we have 
               \begin{equation}\label{thm 10, eq 4}
               	\left( f \left(   M \left( x_{n_i},  x_{n_{i+1}},   {\epsilon}  \right)  \right) \right)^\alpha 
               	\geq  f \left(   M \left( x_{n_i}, x, \frac{\epsilon}{2} \right) \star  M \left(  x_{n_{i+1}}, x,  \frac{\epsilon}{2} \right)  \right). 
               \end{equation}

               This yields 
               \begin{align*}
               	&  (f ( 1-r))^\alpha \geq    f \left(   M \left( x_{n_i}, x, \frac{\epsilon}{2} \right) \star  M \left(  x_{n_{i+1}}, x,  \frac{\epsilon}{2} \right)  \right)  ~  ~~ \forall ~ i \geq N ~~~~ (\text{Using} ~ (\ref{thm 10, eq 2,5}) ~ \& ~ (\ref{thm 10, eq 4}) )\\
               	 \implies & \left(  f ( 1-r) \right)^\alpha >   ( 1-\epsilon)  \quad ( \text{using} ~ (\ref{th 10,eq 1}) ~ \& ~ (\ref{thm 10, eq 3.5})) \\
                 \implies &  \left(  f ( 1-r) \right)^\alpha \geq 1  \quad (\text{since} ~  \epsilon >0 ~  \text{is  arbitrarily chosen}).
               \end{align*}
               
            Hence  using $(\mathscr{F}2 ))$,  we have  
            %%%%%%%%%%%%
            \begin{align*}
                  f (1-r) = 1 \implies  1-r = 1  \implies r = 0.
            \end{align*}
            
    This is a contradiction to our assumption. Hence this case is absurd.    
        This completes the proof. 
         	\end{proof}

% \begin{ppn}\label{ppn 7}
 %	In a fuzzy $\mathcal{F}$-metric space $(X, M, \textit{f}, \alpha,\star)$, for any subset $A \subset X $, 
 %	$x \in$ ${\bar{A}}$, $ 1>r>0 $  implies   $B_\mathcal{F}(x,r,t)$ $\cap A \ne$ $\phi$ for any $ t> 0$.
% \end{ppn} 	
%\begin{proof}
 	%	Suppose $x \in$ ${\bar{A}}$ and $ 1>r>0 $ be  fixed.\\Then $\exists$ a sequence $\{x_n\}$ $\subseteq A$ such that $\lim_{n \to \infty }$ $M(x_n, x,t) = 1$ for any $ t>0 $
 		%\begin{align*}
 		%	&\implies \exists N \in \mathbb{N} ~ such ~ that ~ M(x_n, x,t) > 1-r  ~\forall~ n \geq N\\
 		%	&\implies x_n \in B_\mathcal{F}(x,r,t)~ \forall ~n \geq N\\
 		%	&\implies B_\mathcal{F}(x,r,t) \cap A \ne \phi .
 		%\end{align*}
 	%\end{proof}

  \section{Fixed point theorem  with application to a satellite web coupling problem }
  
  In this section, we present a fixed point theorem in  fuzzy $ \mathscr{F}$- metric space setting that provides a useful criterion for existence of unique fixed point of self-mappings in fuzzy metric spaces that generalizes classical results by incorporating the notion of fuzzy contractions. 
  
  \begin{dfn}
  Let $ ( X, M,  f, \alpha, \star ) $ be a fuzzy $ \mathscr{F} $-metric space. A mapping $ T : X \to X $   is said to be a    fuzzy $ \psi$-contraction mapping with respect to the  function $ \psi \in \Psi $ if the following implication   holds:
  $$ 0 < M ( x, y, t ) <  1 \implies M ( Tx, Ty, t ) \geq \psi ( M ( x, y, t ) ) $$
  for all $ x, y \in X $ and $t> 0 $. 
  \end{dfn}

  \begin{thm} \label{main thm}
  	Let $ ( X, M, f, \alpha, \star  ) $ be an $ \mathscr{F} $-complete fuzzy $ \mathscr{F} $-metric space and   $ T  $     be a    fuzzy $ \psi$-contraction mapping with respect to some  function $ \psi \in \Psi $ such that   there exists $ y \in X $ satisfying  $ 0 < M ( y, T y, t ) <  1 ~ \forall t > 0 $. Then $ T $ has   { at most } one fixed point in $X$. 
  \end{thm}
 \begin{proof}
	Let $ x \in X $ be such that $ 0 < M ( x, T x, t ) <  1 ~ \forall t > 0 $ and define $x_n = T^n x, ~ n \in \mathbb{N} \cup \{ 0 \} $. If $ T ( x_k ) = x_k $ for some $ k \in \mathbb{N} \cup \{ 0 \} $ then the proof is done.  So suppose $ T ( x_r ) \ne  x_r $ for all $ r \in \mathbb{N} \cup \{ 0 \} $. 
    
	Since $T$ is a fuzzy $ \psi$-contraction mapping then for the sequence $ \{ x_n \} $, we have 
	$$ M ( x_{n+2}, x_{n+1}, t ) \geq \psi \left( M ( x_{n+1}, x_{n}, t )\right) ~~  \forall n \in \mathbb{N} \cup \{ 0 \} ~ \& ~ t > 0.  $$ 

	Repeating the above relation, for any $ m \geq 2 $ we get
	$$ M ( x_{m+1}, x_m, t ) \geq  \psi ( M ( x_{m }, x_{m-1}, t ) ) \geq \psi^2 ( M ( x_{m - 1}, x_{m- 2}, t ) )  \geq \cdots \geq \psi^m  ( M ( x_{ 1}, x_{0}, t ) ) ~~ \forall t> 0.  $$ 

	Taking limit as $ m \to \infty $ and using Lemma \ref{prop of si 2} in the above relation, we get 
	\begin{align}\label{main thm 1}
		&  \underset{m \to \infty  }{\lim } M ( x_{m+1}, x_m, t ) \geq \underset{m \to \infty  }{\lim } \psi^m  ( M ( x_{ 1}, x_{0}, t ) )  = 1 ~~ ~~ \forall t> 0         \nonumber \\
		\implies & \underset{m \to \infty  }{\lim } M ( x_{m+1}, x_m, t ) = 1  ~~ ~~ \forall t> 0.
	\end{align} 

	Proceeding in this way, we obtain 
	\begin{equation} \label{main thm 4}
		 \underset{m \to \infty  }{\lim } M ( x_{m+i}, x_{m+i -1}, t ) = 1  ~~ ~~ \forall t> 0 ~ \& ~ i \geq 1.
	\end{equation} 
	Next we show that $ \{ x_r \}$ is a Cauchy sequence in $ ( X, M, f, \alpha, \star  ) $. 
    
	Now using the inequality $ ( \mathscr{F}M4)$, for any $ m  \in \mathbb{N} $ and $ p = 1, 2, 3, \ldots $,  we have
	\begin{equation} \label{main thm 5} 
		(f ( M ( x_{ m+ p }, x_{m}, t ) ))^\alpha \geq f ( M ( x_{m+p}, x_{m+p-1}, \frac{ t }{  p } ) \star \cdots \star M ( x_{m+ 2}, x_{m + 1}, \frac{ t }{ p} ) \star M ( x_{m+1}, x_m, \frac{ t }{ p } ) )    
	\end{equation}

	Applying the relation (\ref{main thm 4}) on (\ref{main thm 5}), we get  
	 $$ \underset{m \to \infty  }{\lim }\{ M ( x_{m+p}, x_{m+p-1}, \frac{ t }{ p } ) \star \cdots    \star M ( x_{m+1}, x_m, \frac{ t }{ p} ) \} = 1,  ~ \forall t> 0 ~ \& ~ p = 1,2, 3, \ldots. $$

	 This implies 
	\begin{align*}
	  & \underset{m \to \infty  }{\lim } f ( M ( x_{m+p}, x_{m+p-1}, \frac{ t }{ p } ) \star \cdots   \star M ( x_{m+1}, x_m, \frac{ t }{ p} )  ) = 1,  ~ \forall t> 0 ~ ~ (\text{by } ~ (\mathscr{F}2)) \\
	or ~~ & 	\underset{m \to \infty  }{\lim } ~ (f ( M ( x_{ m+ p }, x_{m}, t ) ))^\alpha \geq 1 ~~~ \forall t> 0 ~~ (\text{by} ~( \ref{main thm 5})) \\
	or ~~ & \underset{m \to \infty  }{\lim } ~ (f ( M ( x_{ m+ p}, x_{m }, t ) ))^\alpha = 1 ~~~ \forall t> 0 \\
	or ~~ & \underset{m \to \infty  }{\lim } ~ (f ( M ( x_{ m+p }, x_{m}, t ) ))  = 1 ~~~ \forall t> 0 \\
	or ~~ & \underset{m \to \infty  }{\lim } ~  M ( x_{ m + p }, x_{m}, t )  = 1 ~~~ \forall t> 0. 
\end{align*} 

This proves that $ \{ x_m \}$ is a Cauchy sequence in $ ( X, M, f, \alpha, \star  ) $.  Since $X$ is $ \mathscr{F}$-complete, so $ \{ x_m \}$ converges to some $ u \in X $.  
Therefore,
$$ \underset{m \to \infty  }{\lim } ~  M ( x_{ m  }, u, t )  = 1 ~~~ \forall t> 0.  $$

Next we claim that $ u   $ is a fixed point of $T$. If possible suppose there exists $ t_0 > 0$ such that 
\begin{equation}\label{main thm 2}
	M ( u, Tu, t_0 ) < 1. 
\end{equation}

Then we have,
\begin{equation}\label{main thm 3}
	 (	f ( M ( u, Tu, t_0 ) ))^\alpha \geq f ( M ( u, x_m, \frac{t_0}{2} ) \star M ( x_m, Tu, \frac{t_0}{2} )) ~~~\forall m \in \mathbb{N}.
\end{equation}

Again,
\begin{align*}
  &   M ( x_m, Tu, \frac{t_0}{2} ) \geq \psi ( M ( x_{m - 1 }, u, \frac{t_0}{2} ) ) ~~~\forall m\in \mathbb{N} \\ 
\implies & \underset{m \to \infty  }{\lim } ~ M ( x_m, Tu, \frac{t_0}{2} ) \geq \underset{m \to \infty  }{\lim }  ~ \psi ( M ( x_{m - 1 }, u, \frac{t_0}{2} ) ) = 1 ~~~(\text{using  Lemma} ~ \ref{prop of si 2}) \\
\implies & \underset{m \to \infty  }{\lim } ~ M ( x_m, Tu, \frac{t_0}{2} )  = 1. 
\end{align*}

 Therefore,  
 \begin{align*}
 & \underset{m \to \infty  }{\lim } ~ M ( u, x_m, \frac{t_0}{2} ) \star M ( x_m, Tu, \frac{t_0}{2} )  = 1 \\
\implies & \underset{m \to \infty  }{\lim } ~ f ( M ( u, x_m, \frac{t_0}{2} ) \star M ( x_m, Tu, \frac{t_0}{2} )) = 1 ~~~ (\text{using}~ (\mathscr{F}2))  \\
\implies & \underset{m \to \infty  }{\lim } ~  (	f ( M ( u, Tu, t_0 ) ))^\alpha \geq 1 ~~~ (\text{using the relation }~  (\ref{main thm 3}))   \\
\implies &   M ( u, Tu, t_0 ) = 1   ~~~ (\text{using}~ (\mathscr{F}2)).
\end{align*}

Thus we arrived at a contradiction to the relation (\ref{main thm 2}). Therefore $ u $ is a fixed point of $T$. 

To prove the uniqueness of the fixed point, suppose there exists an element $ v \in X $ with $ u \neq v $ such that $T v = v $. Since $ v \neq u $, so there exists $ s > 0 $ such that $ M ( u, v, s ) < 1 $. In that case, we get 
$$ M ( u, v, s ) = M ( T u, T v, s ) \geq \psi ( M ( u, v, s ) ) > M ( u, v, s ), $$
which is a contradiction. So, $T$ has unique fixed point in $X$. 
 \end{proof}

We present the following example in the support of the  above theorem.

 \begin{eg}
 	Consider the Example \ref{ex from 1st paper} of fuzzy $ \mathscr{F} $-metric space. Here we take  
 	 $ X = [ - 5, 5 ] $. 
 	At first we show that $ ( X, M, f, \alpha, \star )$ is complete fuzzy $ \mathscr{F} $-metric space.  For, consider a Cauchy  sequence $ \{ x_n \} $ in $ ( X, M , f, \alpha, \star ) $. 
 	Then 
 	\begin{align*}
 		& \underset{ m, n \to  \infty }{\lim} ~  M ( x_n, x_m, t) = \underset{ m, n \to  \infty }{\lim} ~ \left( \frac{t}{t+1} \right)^{|x_n-x_m|^2} = 1,   ~  ~ ~ \forall t \in ( 0, \infty ). \\
 		\implies & \underset{ m, n \to \infty }{\lim} ~ {|x_n-x_m|}^2 = 0.  \\
 	\implies 	&  \underset{ m, n \to \infty }{\lim} ~ | x_n-x_m| = 0.
 	\end{align*} 
 	Therefore, $ \left\{ x_n \right\} $ is a Cauchy sequence in $ X $ with respect to usual metric .  
 	Since $ X $ is complete with respect to usual metric, so there exists $ x \in X $ such that $ \underset{ n \to \infty } \lim | x_n - x |  = 0 $ which implies
 	  $$ 	\underset{ n \to \infty } \lim M (x_n, x, t)    = \underset{ n \to  \infty } \lim  \left(   \frac{t}{1+t} \right)^{ | x_n - x|^2 }   = 1,  ~ ~ ~  \forall  t \in ( 0, \infty ). $$
 This implies  $ x_n \to x ~\text{as } n \to \infty $ in $ ( X, M, f, \alpha, \star )$.   
 	Hence $ ( X, M, \textit{f}, \alpha, \star) $ is a $ \mathscr{F} $-complete fuzzy $ \mathscr{F} $-metric space.
    
 	Next we consider two mappings   by $ T (x) = \frac{x}{6} ~ \forall x \in X $ and  $ \psi ( t ) = \sqrt{t} ~ \forall t \in [ 0, 1 ] $.
 	Then we have
 	 $$ M ( Tx, Ty, t) = M ( \frac{x}{6}, \frac{y}{6}, t ) = \left( \frac{t}{1+t} \right)^{ \frac{|x-y|^2}{36}} ~~ \text{and} ~~ 
  \psi \left(  M ( x, y, t ) \right) = \sqrt{ M ( x, y, t ) } = \left(  \frac{t}{1+t} \right)^ { \frac{ | x-y|^2}{2} } $$
   implies
   $$ M ( Tx, Ty, t) \geq \sqrt{ M ( x, y, t )} = \psi \left(  M ( x, y, t ) \right), ~~\text{since} ~ \frac{t}{1+t} < 1. $$
 	Hence, $ T $ satisfies the fuzzy $ \psi $-contraction condition.
 	So by above  theorem  $ T $ has a unique fixed point in $ X $ which is $ x = 0 $.   
    
    The geometric illustration of this example is provided in Figure 1 which shows the variation between $ M ( Tx, Ty, t) $ and $  \psi \left(  M ( x, y, t ) \right) $.  In Figure 1$(a) $,   the variation is due to particular fixed value of $ t $ and  Figure 1$(b) $ shows the variation for fixed values of $ | x - y |, x, y\in X $. This visualization highlights the interplay between the mapping $T$   and the function  $ \psi $ that demonstrates  how the inequality $ M ( Tx, Ty, t) \geq  \psi \left(  M ( x, y, t ) \right) $  is sustained under the given conditions. The figure is not only for supporting  the theoretical framework of the contraction condition but also it provides an intuitive understanding of the relationship between the parameters  and  makes the main result  more comprehensible and accessible.

 	\begin{figure}
		[ht!]
		 \centering
		 (a)  \includegraphics[width=6cm,height=6cm]{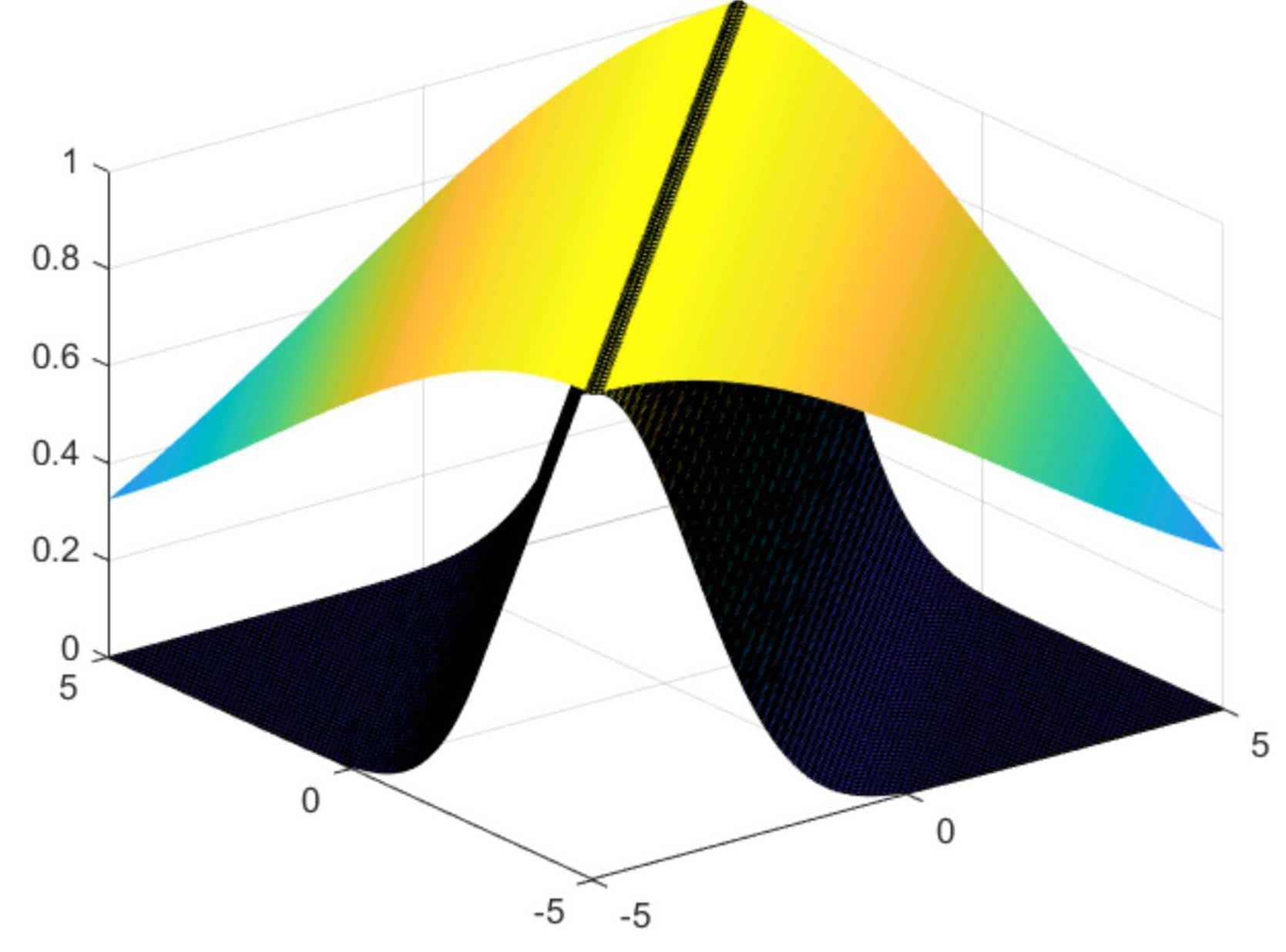}
		%\caption{AAA}
		%
		 (b) \includegraphics[width=6cm,height=6cm]{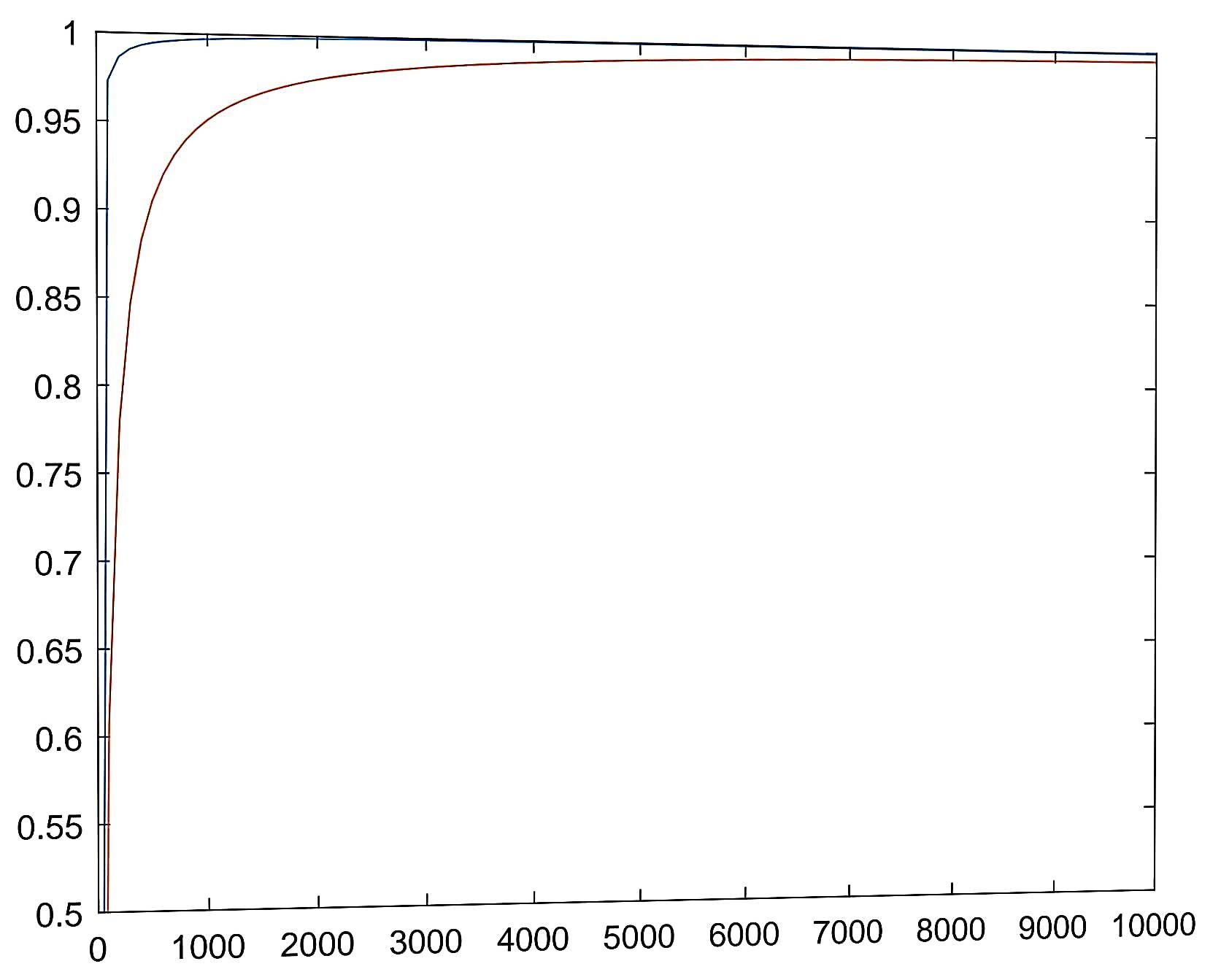}
		\caption{ (a) $ M ( Tx, Ty, t ) $ vs $ \psi ( M (x, y, t )) $ at $ t = 2 $~~~(b) $ M ( Tx, Ty, t ) $ vs $ \psi ( M (x, y, t )) $ for  $ | x -y | = 10 $}
	\end{figure}
  \end{eg} 
  
In the following example we show that the conditions of the theorem are sufficient but  not necessary. 
%	To justify the $ \psi $-contraction condition is not necessary for the above fixed point theorem.  We consider the following example 
            %%%%%%%%%%%
            %%%% 
            
\begin{eg}
Consider the fuzzy $ \mathscr{F} $-metric $M$ of Example \ref{ex from 1st paper} over $ X = \{ 0, 2, 4, \cdots \} $ and define a self mapping $ T : X \to X $   by $ T ( x) = 10x, ~ \forall x \in X $. 

            First we show that $ X $ is complete in $ \mathbb{ R}$ with respect to usual metric . For, let $ \{ x_n \} $ be a Cauchy sequence in $ (X, M, f, \alpha, \star ) $. Then
            $$ \underset{ m, n \to \infty } \lim { M ( x_n, x_m, t) } = 1 ~~ 
            i.e. ~  \underset{ m, n \to \infty } \lim { | x_n - x_m |^2 } = 0 ~ ~ 
             i.e. ~   \underset{ m, n \to \infty } \lim { | x_n - x_m | } = 0. $$   Thus, $ \{  x_n \}$ is Cauchy sequence $ X $ in with respect to usual metric.
            Since $ X $ is closed in $ \mathbb{R} $ under the usual metric.     So, $ X $ is complete thus  $ \exists  x \in  X $ such that 
            $$  \underset{ n \to \infty } \lim { | x_n - x | } = 0 ~     \implies  \underset{ n \to \infty } \lim { | x_n - x |^2 }  = 0 $$
            and hence   $ \underset{ n \to \infty } \lim { M ( x_n, x, t ) } = \underset{ n \to \infty } \lim \left( { \frac{t}{t+1} }\right)^{ | x_n - x |^2 } = \left({\frac{t}{t+1}}\right)^{ \underset{ n \to \infty } \lim { |x_n-x|^2 } } = 1 $.  
            
           Therefore, $ ( X, M, \textit{f}, \alpha, \star ) $ is $ \mathscr{F} $-complete fuzzy $ \mathscr{F} $-metric space. 
           
           Next consider the mapping $ \psi ( t ) = \sqrt{t}, ~ \forall t \in [ 0, 1 ] $. 
              Now observe that, 
             \begin{align*}
             	M ( Tx, Ty, t ) & = \left(  \frac{ t }{ t + 1} \right)^{ | Tx - Ty |^2 } 
                 	= \left(  \frac{ t }{ t + 1 } \right)^{ |10x - 10y|^2 } 
              	= \left( \frac{ t }{ t + 1 } \right)^{ 100 | x- y |^2 }
             \end{align*}
            $$   and  ~~  \sqrt{ M ( x, y, t ) } = \left(   \frac{ t }{ t + 1 } \right)^{ \frac{ | x-y |^2}{2} }. $$
           Since $ | x -y |^2 \geq 4 $ for $ x \ne y $, then $ \frac{| x-y|^2}{ 2 } \leq 100 | x -y |^2 $. 
           So we can write
           \begin{align*}
           	& \left(  \frac{ t }{ t + 1} \right)^{ 100 | x-y|^2 } \leq \left( \frac{ t }{ t + 1 } \right)^{ \frac{ | x-y|^2 }{ 2 } } \\
           	\implies & M ( Tx, Ty, t ) \leq \sqrt{ M ( x, y, t ) }  = \psi ( M (x,y,t) ). 
           \end{align*}

          This shows that $ T $ does not satisfy the contraction condition but it has a  fixed point in $ X $ which is $`0$'. 
\end{eg}

\subsection{Application to a satellite web coupling problem:}  

Inspired by the applications of fixed point techniques to a variety of real-world problems, we employ Theorem \ref{main thm} to address a satellite web coupling boundary value problem \cite{new-8.1.25}. A satellite web coupling can be conceptualized as a thin sheet linking two cylindrical satellites. The analysis of radiation emanating from the web coupling between these satellites gives rise to the following nonlinear boundary value problem:
    \begin{equation}\label{ satellite 1 }
              	- \frac{ d^2 \omega }{ d t^2 } = \mu \omega^4, ~ 0 < t < 1, ~  \omega ( 0 ) = \omega ( 1 ) = 0 
    \end{equation}
    where  $\omega(t) $  represents the temperature of radiation at any point $ t \in [ 0, 1 ] $, $ \mu = \frac{ 2 a l^2 K^3 }{ \zeta h } > 0 $ is a dimensionless positive constant, $ K $ denotes the constant absolute temperature of both satellites, while heat is radiated from the surface of the web at $0 $ absolute temperature, $l$ is the distance between the two satellites, $a$ is a positive constant characterizing the radiation properties of the web's surface, the factor $2$ accounts for radiation from both the top and bottom surfaces,  $ \zeta $  is the thermal conductivity, and $h$ is the thickness of the web.
    
Next we recall   the corresponding integral equation: 
       $$   \omega ( t ) = 1 - \mu \int_{ 0 }^{ 1 } G ( t, \xi ) \omega^4 ( \xi ) d\xi  $$ 
where $ G ( t, \xi ) $ is the Green function, defined by
    $$ G ( \alpha, \xi ) = \begin{cases}
         	\alpha ( 1 - \xi ), \quad 0 < \alpha < \xi \\
         	\xi ( 1 - \alpha ), \quad \xi < \alpha < 1. 
         	\end{cases} $$

Let  $  {X} = C[ 0, 1 ] $, the class of all real valued continuous functions defined on $ [ 0, 1 ] $ and define a function $ \mathcal{M} :  {X} \times  l{X} \times ( 0, \infty )  \to [ 0, 1 ] $ by 
$$  {M} ( \xi, \eta, T ) = \left( \frac{ T }{ T + 1 } \right)^{ \underset{s \in [ 0, 1 ] }{\sup}| \xi ( s ) - \eta ( s ) |^2}$$
for all $ \xi, \eta \in  {X} $ and $ T > 0 $. Then following the  Example \ref{ex from 1st paper}, we can show that $ {M}$ is a fuzzy $ \mathscr{ F } $-metric on $X$ with $ f ( x )= x^2, ~ \forall x \in [ 0, 1 ]$, $ \alpha = \frac{ 1 }{ 2 } $  and $ \star = \text{product} $. Moreover, $(  {X},  {M}, f, \alpha, \star ) $ is a $ \mathscr{ F } $-complete fuzzy $ \mathscr{ F } $-metric space.\\

\begin{thm}
Consider the $ \mathscr{ F } $-complete fuzzy metric space $ ( X, M, f, \alpha, \star ) $defined above. Suppose that the boundary value problem satisfies the following condition: 
$$ \underset{ \xi \in [0,1] } \sup | \left( \omega^2 ( \xi ) + v^2 ( \xi ) \right) \left( \omega ( \xi ) + v ( \xi ) \right)| \leq \frac{ k }{ \mu } ~ ~\quad \text{where} ~ k \in ( 0, 4 ) $$ 
 Then the satellite web coupling boundary value problem (\ref{ satellite 1 }) has a unique solution.
 \end{thm}
\begin{proof}
      We define a self mapping $ A : X \to X $ by 
               	 $$  A ( \omega ( t ) ) =  1 - \mu \int_{ 0 }^{ 1 } G ( t, \xi ) \omega^4 ( \xi ) d\xi, ~ \xi \in [ 0, 1 ] $$
    Clearly, a solution to the satellite web coupling problem (\ref{ satellite 1 }) corresponds to a fixed point of the self-mapping $ A $. 
    
    Now for all $ \omega, v \in X $ and $ t \in [ 0, 1 ] $,
               	 \begin{align*}
               	 	& | A \omega ( t ) - A v ( t ) |^2 \\
                    = & ~ \mu^2 \bigg| \int_{ 0 }^{ 1 } ( \omega^4 ( \xi ) - v^4 ( \xi ) ) ~ G ( t, \xi ) d\xi \bigg|^2 \\
               	 = & ~   \mu^2 \bigg| \int_{ 0 }^{ 1 } \left\{ \left( \omega^2 ( \xi ) + v^2 ( \xi ) \right) \left(   \omega ( \xi ) + v ( \xi )  \right) \left( \omega ( \xi ) - v ( \xi ) \right)  \right\}  ~    G ( t, \xi ) d\xi \bigg|^2 \\
               	 	 \leq & ~ \mu^2 ~   \underset{ s \in [ 0, 1 ] }{\sup} \big| ( \omega^2 ( s ) + v^2 ( s ) )  (  \omega ( s ) + v ( s ) ) \big|^2 ~  \bigg| \int_{ 0 }^{ 1 }   ( \omega ( \xi ) - v ( \xi ) )  ~ G ( t, \xi ) d\xi  \bigg|^2 \\
               	 	 \leq & ~ \mu^2 ~ \frac{ k^2 }{ \mu^2 } \underset{ s \in [ 0, 1] }{ \sup} \big| \omega ( s) - v ( s ) \big|^2 ~ \bigg| \int_{ 0 }^{ 1 } G ( t, \xi ) d\xi \bigg|^2  \\
               	 	\leq &  ~ k^2 ~ \underset{ s \in [ 0, 1] }{ \sup} \big| \omega ( s ) - v ( s ) \big|^2 ~ \underset{ t \in [ 0, 1] }{ \sup}  \bigg| \int_{ 0 }^{ 1 } G ( t, \xi ) d\xi   \bigg|^2 \\
               	 	 = & ~ \frac{ k^2 }{ 16 } \cdot \underset{ s \in [ 0, 1] }{ \sup} \big| \omega (  s ) - v (s) \big|^2 \\
               	 	= &  ~ \frac{ 1 }{ \beta } \cdot \underset{ s \in [ 0, 1] }{\sup} | \omega ( s ) - v ( s ) |^2  \quad  \text{where} ~ \frac{ 1 }{ \beta } = \frac{ k^2 }{ 16 }   \in ( 0, 1 ), ~ \beta > 0.
               	 \end{align*}
              This implies  
              \begin{align*}
                  & \underset{ t \in [ 0, 1 ] }{ \sup} | A \omega ( t ) - A v ( t ) |^2 \leq \frac{ 1 }{ \beta } \cdot \underset{ s \in [ 0, 1 ] }{ \sup} | \omega  ( s ) - v( s ) |^2 \\
                  or & ~ \left(  \frac{ T }{ T + 1 } \right)^{ \underset{ t \in [ 0, 1 ]} \sup | A \omega ( t ) - A v ( t ) |^2 } \geq \left( \frac{ T }{ T +1 } \right)^{ \frac{ 1 }{ \beta } \cdot \underset{ t \in [ 0, 1 ] } \sup  | \omega ( t ) - v ( t ) |^2 } ~~ \forall T > 0.
              \end{align*}

                 Now  consider the mapping $ \psi ( t ) = t^{ \frac{ 1 }{ \beta } } ~~ \forall t \in [ 0, 1] $ where $ \frac{ 1 }{ \beta } = \frac{ k^2 }{ 16 } \in ( 0, 1 ) $. 
                 
                Then the above inequality gives 
                \begin{align*}
                	 &  M ( A \omega, A v, T ) \geq \left(  M ( \omega, v, T )  \right)^{ \frac{ 1 }{ \beta }} ~~ \forall T > 0  \\
                	or ~ ~ &   M ( A \omega, A v, T ) \geq \psi \left( M ( \omega, v, T) \right)  ~~ \forall T > 0.
               \end{align*}	
           
             Thus, the mapping $ A $ fulfill the conditions of the Theorem \ref{main thm} and therefore $ A $ has a unique fixed point in $ X $.  
             Consequently, the boundary value problem (\ref{ satellite 1 }) has a solution in $X$.  
               \end{proof}

  \section*{Conclusion.}
  In this article, we have investigated several fundamental properties  related to the compactness and totally boundedness  of   fuzzy $\mathscr{F}$-metric spaces.     Within this framework, we prove a   fixed point theorem which not only enriches the existing literature but also demonstrates practical applicability through its relevance to the satellite web coupling problem. To illustrate the sufficiency of the theorem’s conditions, we present concrete examples.
  A geometric illustration of the example is also provided,  which shows the variation between $ M ( Tx, Ty, t) $ and $  \psi \left(  M ( x, y, t ) \right) $ and highlights shows the   interplay between the mapping $T$   and the function  $ \psi $. The figure  stands to support  the theoretical foundation of the contraction condition and offers an intuitive understanding of the relationships among the parameters involved.
  
 These results lay a foundation for further research  with the possibility of extensions  to study more generalized spaces or exploring additional applications in dynamic systems, optimization  and computational mathematics.

 ~\\

%  \paragraph{Author Contributions.} All Authors (DB, AD and TB) have contributed as follows: conceptualization, investigation, methodology, writing–original draft, DB; 
  %
%  conceptualization,  methodology,  writing–review \& editing,  AD;
  %  
 % conceptualization,  validation, writing–review \& editing , T. Bag.  All authors have read and approved the published version of the manuscript.

  % ----------------------------------------------------------------
  
  \paragraph{Conflicts of Interest.} The authors declare no conflicts of interest. \\
  
  % ----------------------------------------------------------------
     
  \noindent{\bf Acknowledgment.} 
  The author DB is thankful to University Grant Commission (UGC), New Delhi, India for awarding her junior research fellowship [Ref. No 231610065558 (CSIR-UGC NET DECEMBER-2022/JUNE-2023)].
  The   authors also   are thankful to the Department of Mathematics, Siksha-Bhavana, Visva-Bharati, India.

\end{document}